\documentclass[twocolumn]{autart}    

\usepackage{graphicx} 
\usepackage[dvipsnames]{xcolor}

\usepackage{amsmath,amssymb,amsfonts}
\usepackage{algorithm}
\usepackage{algpseudocode}
\usepackage{array}
\usepackage{bbm}
\usepackage{balance}
\usepackage{booktabs}
\usepackage{cases}
\usepackage{cite}
\usepackage{dsfont}
\usepackage[official]{eurosym}
\usepackage{epstopdf}
\usepackage[acronym]{glossaries}
\usepackage{mathrsfs}
\usepackage{mathtools}
\usepackage{multirow}
\usepackage{makecell}
\usepackage{textcomp}
\usepackage{soul}
\usepackage{subfigure}
\usepackage{url}
\usepackage{siunitx}
\usepackage{changes}

\newtheorem{clr}{Corollary}
\newtheorem{theorem}{Theorem}

\newtheorem{remark}{Remark}

\newenvironment{proof}{{\noindent\it Proof.}\quad}{\hfill $\square$\par}

\begin{document}

\renewcommand{\algorithmicrequire}{\textbf{Input:}}
\renewcommand{\algorithmicensure}{\textbf{Output:}}

\begin{frontmatter}

\title{Data-Driven Output Prediction and Control of Stochastic Systems: An Innovation-Based Approach\thanksref{footnoteinfo}} 

\thanks[footnoteinfo]{This work was supported by National Natural Science Foundation of China under Grants 62373211, 62325305 and 62327807. This paper was not presented at any IFAC meeting. Corresponding author C. Shang. Tel. +86-10-62782459. Fax +86-10-62773789.}
\author[Tsinghua]{Yibo Wang}\ead{wyb21@mails.tsinghua.edu.cn},
\author[Tsinghua]{Keyou You}\ead{youky@tsinghua.edu.cn},
\author[Tsinghua]{Dexian Huang}\ead{huangdx@tsinghua.edu.cn},
\author[Tsinghua]{Chao Shang}\ead{c-shang@tsinghua.edu.cn},
\address[Tsinghua]{Department of Automation, Beijing National Research Center for Information Science and Technology, Tsinghua University, Beijing 100084, China}


\begin{abstract}                          
Recent years have witnessed a booming interest in data-driven control of dynamical systems. However, the implicit data-driven output predictors are vulnerable to uncertainty such as process disturbance and measurement noise, causing unreliable predictions and unexpected control actions. In this brief, we put forward a new data-driven approach to output prediction of stochastic linear time-invariant (LTI) systems. By utilizing the innovation form, the uncertainty in stochastic LTI systems is recast as innovations that can be readily estimated from input-output data without knowing system matrices. In this way, by applying the fundamental lemma to the innovation form, we propose a new innovation-based data-driven output predictor (OP) of stochastic LTI systems, which bypasses the need for identifying state-space matrices explicitly and building a state estimator. The boundedness of the second moment of prediction errors in closed-loop is established under mild conditions. The proposed data-driven OP can be integrated into optimal control design for better performance. Numerical simulations demonstrate the outperformance of the proposed innovation-based methods in output prediction and control design over existing formulations.
\end{abstract}

\end{frontmatter}

\section{Introduction}
In modern control applications, the explosive complexities pose a challenge to the traditional identification-for-control scheme, since accurate modeling is costly and heavily relies on engineering expertise. It is thus attractive to design controllers directly from raw data with the identification step bypassed, considering the ever-growing availability of data. From this new perspective, a paradigm shift has been driven towards data-driven end-to-end solutions to output prediction and control design \cite{hou2013model}. A mainstream of research builds upon a notable result rooted in the behavioral theory, later known as the \textit{Willems' fundamental lemma} \cite{willems2005note}. It dictates that in the deterministic case, all trajectories from a linear system are expressible as linear combinations of known trajectories provided the inputs are persistently exciting of a sufficiently high order. This trait leads to the well-known Data-enabled Predictive Control (DeePC) \cite{coulson2019data}, where the raw data-based constraints are used in lieu of parametric models for the description of system dynamics and control design. Thanks to these merits, DeePC and its variants have found miscellaneous applications in power systems \cite{huang2021decentralized,bilgic2022data}, motion control \cite{elokda2021data,carlet2022data}, and smart buildings \cite{chinde2022data}.

Admittedly, prediction accuracy is pivotal for the success of predictive control design. However, the philosophy of substituting state-space models with data-based trajectory representations, in the spirit of the fundamental lemma, is only valid in an ideal deterministic context. In real-life applications, data-driven prediction becomes highly vulnerable to noise corruption and thus the fidelity of DeePC is greatly challenged. A series of remedies have been put forward, making heavy use of regularization \cite{coulson2019data,dorfler2022bridging}. Recently, a new route towards data-driven control of stochastic systems has emerged, by regarding additive uncertainties as exogenous inputs and invoking the fundamental lemma; see e.g. \cite{huang2021decentralized,pan2021stochastic}. 

In this brief, we consider data-driven prediction and control of stochastic linear time-invariant (LTI) systems, where states, additive process disturbance, and measurement noise are unmeasurable. The conventional model-based control relies on the design of Kalman filter (KF), which requires knowing state-space matrices and noise statistics. In a data-driven setting, we recapture additive uncertainties in terms of innovations, i.e., one-step prediction errors in the steady-state KF (SSKF), by converting the stochastic system to its innovation form, where the innovations can be readily estimated from input-output data without knowing state-space matrices, by simply fitting a non-parametric vector auto-regressive with exogenous input (VARX) model. This enables us to apply the fundamental lemma to the innovation form and then develop a new data-driven output predictor (OP) of the stochastic LTI system, which ``implicitly" identifies state-space matrices and performs state estimation. We establish the condition when our data-driven OP and a model-based SSKF perform equally, which then inspires an easy closed-loop implementation of data-driven OP. More importantly, the boundedness of the second moment of prediction errors is established, provided that the input is bounded, the output has a bounded second moment, and a matrix constructed from attainable data is Schur stable. The boundedness offers a practical guideline to verify the ``validity" of innovation estimates as well as the data-driven OP. Finally, the usage of the data-driven OP in control design yields a new innovation-based DeePC (Inno-DeePC) scheme. Numerical examples showcase the performance improvement of the proposed data-driven OP and control over known formulations.

\textbf{Notation:} We denote by $\mathbb{Z}$ ($\mathbb{Z}^+$) the set of (positive) integers. The identity matrix of size $s$ is $I_s\in\mathbb{R}^{s\times s}$. The zero vector of size $s$ and matrix of size $s_1\times s_2$ are denoted by $0_{s}\in\mathbb{R}^s$ and $0_{s_1\times s_2}\in\mathbb{R}^{s_1\times s_2}$. For a matrix $X$, $X^\dagger$ denotes the Moore-Penrose inverse, $\|X\|_F$ is the Frobenius norm, $\|X\|_2$ is the spectral norm, and $X_{[i:j]}$ denotes the submatrix constructed with entries from the $i$-th to the $j$-th row of $X$. Given a sequence $\{ x(i) \}_{i=1}^N\in\mathbb{R}^n$, $x_{[i:j]}$ denotes the restriction of $x$ to the interval $[i,j]$ as $\mathbf{col}(x(i),x(i+1),\cdots,x(j))=[x(i)^\top~x(i+1)^\top~\cdots~x(j)^\top]^\top$.
Similarly, we use $\mathbf{col}(X_1,X_2,\cdots,X_n) = [X_1^\top~X_2^\top~\cdots~X_n^\top]^\top$. A block Hankel matrix of depth $s$ can be constructed from ${x}_{[i:j]}$ via the block Hankel matrix operator $\mathcal{H}_s(x_{[i:j]})$. A sequence $x_{[i:j]}$ is said to be persistently exciting of order $s$ if $\mathcal{H}_s(x_{[i:j]})$ has full row rank. Given an order-$n$ state-space model $(A,B,C,D)$, its lag is denoted as $\ell(A,B,C,D)$, i.e., the smallest integer $l$ such that $\mathbf{col}(C,...,CA^{l-1})$ has rank $n$.

\section{Main Results}

Consider a discrete-time linear time-invariant (LTI) system subject to additive uncertainty:
\begin{equation}
    \label{equation: LTI system}
    \left \{    \begin{aligned}
    x(t+1)&=Ax(t)+Bu(t)+w(t)\\
    y(t)&=Cx(t)+Du(t)+v(t)
    \end{aligned} \right.,
\end{equation}
where ${u}(t)\in\mathbb{R}^{n_u}$, ${y}(t)\in\mathbb{R}^{n_y}$ and ${x}(t)\in\mathbb{R}^{n_x}$ stand for input, output, and unmeasurable state, respectively. $w(t) \in\mathbb{R}^{n_x}$ and $v(t) \in \mathbb{R}^{n_y}$, which are zero-mean white Gaussian noise with covariance matrices $\Sigma_w$ and $\Sigma_v$, denote unmeasurable process and measurement noise that is uncorrelated with $u(t)$ but possibly correlated with $u(t+k),~k>0$. Meanwhile, \eqref{equation: LTI system} is assumed to be minimal, i.e., $(A,C)$ is observable, $(A, B)$ is controllable and $(A,\Sigma_w^{1/2})$ is also controllable. When state measurements are not available, an SSKF can be implemented for state estimation \cite{kamen1999introduction}:
\begin{equation}
\label{equation: steady-state kalman filter}
    \left\{
    \begin{aligned}
        \hat{x}(t+1)&=A\hat{x}(t)+Bu(t)+K[y(t) - \hat{y}(t)]\\
        \hat{y}(t)&=C\hat{x}(t)+Du(t),
    \end{aligned}\right.
\end{equation}
where $K$ is the steady-state Kalman gain rendering all eigenvalues of $\bar{A}\triangleq A-KC$ strictly inside the unit circle. The output predictor is at the heart of predictive control tasks. Note that $\hat{x}(t)=\hat{x}(t|t-1)$ is the state estimate based on available information till time $t-1$. Given $\hat{x}(t)$ from \eqref{equation: steady-state kalman filter} and future inputs, multi-step prediction of state and output trajectories can be performed in a deterministic sense for $k\ge0$:
\begin{equation}
    \label{equation: model-based simulation}
    \left\{ \begin{aligned}
       &\hat{x}(t+k+1|t-1)=A\hat{x}(t+k|t-1)+Bu(t+k)\\
       &\hat{y}(t+k|t-1)=C\hat{x}(t+k|t-1)+Du(t+k),
    \end{aligned} \right.
\end{equation}
where $\hat{x}(t+k|t-1)$ and $\hat{y}(t+k|t-1)$ denote $(k+1)$-step ahead predictions of state and output. More formally, we may refer to \eqref{equation: model-based simulation} as the SSKF-based OP of stochastic systems, which hinges on the knowledge of $(A,B,C,D)$ and noise statistics to derive the steady-state gain $K$. Using $K$ in SSKF \eqref{equation: steady-state kalman filter}, the stochastic LTI system \eqref{equation: LTI system} can be recast into the so-called \textit{innovation form} \cite{huang2008dynamic}:
\begin{equation}
    \label{equation: innovation-form LTI system}
    \left \{ \begin{aligned}
        \hat{x}(t+1)&=A\hat{x}(t)+Bu(t)+Ke(t)\\
        y(t)&=C\hat{x}(t)+Du(t)+e(t),
    \end{aligned}\right.
\end{equation}
where the innovation $e(t)=y(t)-\hat{y}(t)$, defined as the one-step prediction error of SSKF \eqref{equation: steady-state kalman filter}, is a zero-mean white noise process.

\subsection{Data-driven output predictor of stochastic systems}

In this brief, our interest is in predicting the output of \eqref{equation: LTI system} in a \textit{data-driven} fashion, without knowing $(A,B,C,D)$ and noise statistics. Before proceeding, we recall that for system \eqref{equation: LTI system} with $w(t) = 0$ and $v(t) = 0$, a data-driven non-parametric system description can be obtained by invoking the fundamental lemma.

\begin{theorem}[Fundamental Lemma, \cite{willems2005note}]
    Consider the system \eqref{equation: LTI system} with $w(t)=0$ and $v(t)=0$, from which an input-output trajectory $\{u^{\rm d}(i), y^{\rm d}(i)\}_{i=1}^N$ is generated offline. If $N\ge(n_u+1)L+n_x-1$ and $u^{\rm d}_{[1:N]}$ is persistently exciting of order $(L+n_x)$, the following statements hold.
    \begin{enumerate}
        \item[(a)] $\{ u(i), y(i) \}_{i=1}^L$ is a valid input-output data trajectory of system \eqref{equation: LTI system} with $w(t)=0$ and $v(t)=0$ if and only if there exists a vector ${g}\in\mathbb{R}^{N-L+1}$ satisfying $\mathbf{col}(U_{\rm d},Y_{\rm d})g=\mathbf{col}(u_{[1:L]},y_{[1:L]})$,         where ${U}_{\rm d} = \mathcal{H}_L(u^{\rm d}_{[1:N]})$ and ${Y}_{\rm d} = \mathcal{H}_L(y^{\rm d}_{[1:N]})$ are block Hankel matrices of inputs and outputs.
        \item[(b)] Consider the past input-output data $u_{\rm p}(t)= u_{[t-L_{\rm p}:t-1]}$ and $y_{\rm p}(t)= y_{[t-L_{\rm p}:t-1]}$ and the future input $u_{\rm f}(t)= u_{[t:t+L_{\rm f}-1]}$, under the conditions of $L_{\rm p}\ge\ell(A,B,C,D)$ and $L=L_{\rm p}+L_{\rm f}$ where $L_{\rm f}\ge1$ is the future horizon, the future output $y_{\rm f}(t)=y_{[t:t+L_{\rm f}-1]}$ can be uniquely decided by $Y_{\rm f}g(t)=y_{\rm f}(t)$, where $g(t)\in\mathbb{R}^{N-L+1}$ is an additional variable solving the following equations:
        \begin{equation}
            \label{equation: DeePC division}
                \mathbf{col}(U_{\rm p},U_{\rm f},Y_{\rm p})g(t)=
                \mathbf{col}(u_{\rm p}(t),u_{\rm f}(t),y_{\rm p}(t))
        \end{equation}
        with $U_{\rm p} = U_{{\rm d},[1:n_uL_{\rm p}]}\in\mathbb{R}^{n_uL_{\rm p}\times (N-L+1)},~
                U_{\rm f} = U_{{\rm d},[n_uL_{\rm p}+1:n_uL]}\in\mathbb{R}^{n_uL_{\rm f}\times(N-L+1)}$,
        and $Y_{\rm p},Y_{\rm f}$ are constructed in a similar way.
    \end{enumerate}
\label{thm: fundamental lemma}
\end{theorem}

In virtue of the fundamental lemma, a data-driven characterization of innovation form \eqref{equation: innovation-form LTI system} can be obtained by regarding innovations as exogenous inputs. Furthermore, by recasting system \eqref{equation: LTI system} into its innovation form \eqref{equation: innovation-form LTI system}, a data-driven multi-step OP of stochastic system \eqref{equation: LTI system} becomes possible, and its equivalence with SSKF-based OP can be established under certain conditions.

\begin{theorem}
    \label{theorem: extension of fundamental lemma and equivalence to SSKF}
    Consider the finite-dimensional LTI system of innovation form \eqref{equation: innovation-form LTI system}, where an input-output-innovation sequence $\{u^{\rm d}(i),y^{\rm d}(i),e^{\rm d}(i)\}_{i=1}^N$ is generated offline. If $N\ge(n_u+n_y+1)L+n_x-1$ and the extended input signal $\tilde{u}^{\rm d}_{[1:N]}$ defined by $\tilde{u}^{\rm d}(t)=\mathbf{col}(u^{\rm d}(t),e^{\rm d}(t))$ is persistently exciting of order $(L+n_x)$, the following statements hold.
    \begin{enumerate}
        \item[(a)] $\{u(i),y(i),e(i)\}_{i=1}^L$ is a valid trajectory of system \eqref{equation: innovation-form LTI system} if and only if there exists a vector $g\in\mathbb{R}^{N-L+1}$ such that $\mathbf{col}(U_{\rm d},Y_{\rm d},E_{\rm d})g=\mathbf{col}(u_{[1:L]},y_{[1:L]},e_{[1:L]})$,
        where $E_{\rm d}=\mathcal{H}_L(e^{\rm d}_{[1:N]})$.
        \item[(b)] Consider past data $\{u_{\rm p}(t),y_{\rm p}(t),e_{\rm p}(t)\}$ with $e_{\rm p}(t)=e_{[t-L_{\rm p}:t-1]}$ and the future input sequence $u_{\rm f}(t)$, under the condition of $L_{\rm p}\ge\ell(A,B,C,D)$, a data-driven output predictor of \eqref{equation: LTI system} is given by:
        \begin{equation}
            \label{equation: DeePC innovation division yf(t)}
            \hat{y}_{\rm f}(t)=Y_{\rm f}g(t)
        \end{equation}
        with $g(t)$ solving
        \begin{equation}
            \label{equation: DeePC innovation division g(t)}
            \begin{aligned}
            &\mathbf{col}(U_{\rm p},U_{\rm f},Y_{\rm p},E_{\rm p},E_{\rm f})g(t)\\
            &=\mathbf{col}(u_{\rm p}(t),u_{\rm f}(t),y_{\rm p}(t),e_{\rm p}(t),0_{n_yL_{\rm f}}),
            \end{aligned}
        \end{equation}
        where $E_{\rm p}$ and $E_{\rm f}$ are defined by splitting $E_{\rm d}$, akin to $U_{\rm p}$ and $U_{\rm f}$. The data-driven OP \eqref{equation: DeePC innovation division yf(t)} and \eqref{equation: DeePC innovation division g(t)} is equivalent to the SSKF-based OP \eqref{equation: model-based simulation}, where $e_{\rm p}(t)$ in \eqref{equation: DeePC innovation division g(t)} is constituted by past realizations of innovation $\{e(t-L_{\rm p}),\cdots,e(t-1)\}$ from the SSKF \eqref{equation: steady-state kalman filter}. 
    \end{enumerate}
\end{theorem}
\begin{proof}
    The proof can be made by trivially applying Theorem \ref{thm: fundamental lemma} to the system \eqref{equation: innovation-form LTI system} in an innovation form as well as the SSKF-based predictor \eqref{equation: model-based simulation}, which is in a similar spirit to some known extensions of the fundamental lemma to the case with additive disturbance \cite{huang2021decentralized,pan2021stochastic}. Thus, only a sketch is provided here. Due to assumptions behind system \eqref{equation: LTI system}, $(A,[B~ K])$ is controllable. Thus, by viewing $\mathbf{col}(u(t),e(t))$ as extended inputs and invoking Theorem \ref{thm: fundamental lemma}(a), Theorem \ref{theorem: extension of fundamental lemma and equivalence to SSKF}(a) immediately follows. As for Theorem \ref{theorem: extension of fundamental lemma and equivalence to SSKF}(b), the initial state estimate $\hat{x}(t)$ from SSKF \eqref{equation: steady-state kalman filter} is uniquely determined by the past sequence $\{u(i),y(i),e(i)\}_{i=t-L_{\rm p}}^{t-1}$, given $L_{\rm p}\ge\ell(A,B,C,D)$. Thus, the data-driven OP implicitly decides the same initial condition as the SSKF-based OP \eqref{equation: model-based simulation}. Moreover, the multi-step prediction in \eqref{equation: model-based simulation} can be viewed as the output of \eqref{equation: innovation-form LTI system} by setting future innovations $e_{\rm f}(t)$ to zero. As a result, the dynamics of both past and future sections in \eqref{equation: model-based simulation} can be described using the innovation form \eqref{equation: innovation-form LTI system}. Then invoking Theorem \ref{theorem: extension of fundamental lemma and equivalence to SSKF}(a) yields the desired result.
\end{proof}

\begin{remark} When $u^{\rm d}_{[1:N]}$ is persistently exciting of order $L+n_{x}$ as required by Theorem \ref{thm: fundamental lemma}, it is not difficult to satisfy the requirement for persistent excitation of $\Tilde{u}^{\rm d}_{[1:N]}$ in Theorem \ref{theorem: extension of fundamental lemma and equivalence to SSKF}, considering the white noise property of innovation $e^{\rm d}_{[1:N]}$. Because $e^{\rm d}_{[1:N]}$ is a white noise uncorrelated with $u^{\rm d}_{[1:N]}$, the stacked matrix $\mathbf{col}(\mathcal{H}_{L+n_x}(u^{\rm d}_{[1:N]}),\mathcal{H}_{L+n_x}(e^{\rm d}_{[1:N]}))$ has full row rank almost surely.
\end{remark}

Theorem \ref{theorem: extension of fundamental lemma and equivalence to SSKF}(b) describes the condition where such a data-driven predictor is equivalent to a model-based SSKF. However, the proposed data-driven OP hinges on the true values of $e^{\rm d}_{[1:N]}$ coming from a model-based SSKF \eqref{equation: steady-state kalman filter}, so the realization of \eqref{equation: DeePC innovation division yf(t)} and \eqref{equation: DeePC innovation division g(t)} is as hard as exactly knowing matrices $(A,B,C,D,K)$. To enable data-driven implementation of OP \eqref{equation: DeePC innovation division yf(t)} and \eqref{equation: DeePC innovation division g(t)}, we propose to use data-based innovation estimates $\hat{e}^{\rm d}_{[1:N]}$ in place of $e^{\rm d}_{[1:N]}$, since estimates $\hat{e}^{\rm d}_{[1:N]}$ are readily attainable from input-output data by fitting a non-parametric VARX model without knowing $(A,B,C,D,K)$ \cite{chiuso2007role}. The details are deferred to the Appendix.

When $\{E_{\rm p}, E_{\rm f}\}$ in \eqref{equation: DeePC innovation division g(t)} are replaced by their estimates $\{\hat{E}_{\rm p}, \hat{E}_{\rm f}\}$, the equality constraint \eqref{equation: DeePC innovation division g(t)} becomes:
 \begin{equation}
    \label{equation: DeePC innovation division g(t) with innovation estimates}
    \begin{aligned}
    &\mathbf{col}(U_{\rm p},U_{\rm f},Y_{\rm p},\hat{E}_{\rm p},\hat{E}_{\rm f})g(t)\\
    &=\mathbf{col}(u_{\rm p}(t),u_{\rm f}(t),y_{\rm p}(t),\hat{e}_{\rm p}(t),0_{n_yL_{\rm f}}).
    \end{aligned}
\end{equation}
However, in such conditions, Theorem \ref{theorem: extension of fundamental lemma and equivalence to SSKF} does not hold exactly and it is no longer suitable to use arbitrary solution $g(t)$ to \eqref{equation: DeePC innovation division g(t) with innovation estimates} in the proposed data-driven OP, due to inevitable estimation errors in $\hat{E}_{\rm p}$ and $\hat{E}_{\rm f}$. The pseudo-inverse solution
\begin{equation}
    \label{equation: pinv solution with estimated innovation}
    \begin{aligned}
    g_{\rm pinv}(t) = &~\mathbf{col}(U_{\rm p},U_{\rm f},Y_{\rm p},\hat{E}_{\rm p},\hat{E}_{\rm f})^\dagger \\
    &\quad \cdot \mathbf{col}(u_{\rm p}(t),u_{\rm f}(t),y_{\rm p}(t),\hat{e}_{\rm p}(t),0_{n_yL_{\rm f}})
    \end{aligned}
\end{equation}
can be used to alleviate the error effect. Note that in contrast to \eqref{equation: DeePC innovation division g(t)}, $\hat{e}_{\rm p}(t) = [\hat{e}(t-L_{\rm p})^\top~ \cdots ~ \hat{e}(t-1)^\top]^\top$ is used in \eqref{equation: DeePC innovation division g(t) with innovation estimates} and \eqref{equation: pinv solution with estimated innovation} rather than $e_{\rm p}(t) = [e(t-L_{\rm p})^\top~ \cdots ~ e(t-1)^\top]^\top$, because the data-driven OP does not coincide with SSKF and thus a different innovation sequence $\{ \hat{e}(t) \}$ is generated. This will be further clarified in the next subsection. As such, we arrive at a new data-driven OP $\hat{y}_{\rm f}(t)=Y_{\rm f}g_{\rm pinv}(t)$ with $g_{\rm pinv}(t)$ given by \eqref{equation: pinv solution with estimated innovation}, which bypasses the need for both identification of $(A,B,C,D,K)$ and online implementation of state estimator. 

\begin{remark}
    Note that the data-driven approach to innovation estimation, as elaborated in the Appendix, does not identify parameters $(A,B,C,D,K)$ implicitly or explicitly. This is because Markov parameters (MPs) estimated from VARX modeling, i.e., $\hat{\Phi}_u$ and $\hat{\Phi}_y$ solving \eqref{eq: LSR}, do not exactly correspond to a realization of $(A,B,C,D,K)$. Specifically, in closed-loop subspace identification (SID), one has to perform order selection and then estimate $(A,B,C,D,K)$ from $\hat{\Phi}_u$ and $\hat{\Phi}_y$ through least-squares fitting \cite{van2013closed}, which is inevitably subject to errors. By contrast, the proposed data-driven OP $\hat{y}_{\rm f}(t)=Y_{\rm f}g_{\rm pinv}(t)$ with $g_{\rm pinv}(t)$ given by \eqref{equation: pinv solution with estimated innovation} based on innovation estimates bypasses order selection and parameter estimation in SID, thereby circumventing error accumulation.
\end{remark}

\begin{remark}
    In \cite{huang2021decentralized,pan2021stochastic}, various forms of data-driven control have been developed by extending the fundamental lemma to handle additive disturbances. Differently, our focus is placed on the usage of innovation form as a vehicle for devising a practical data-driven OP of \eqref{equation: LTI system}. As to be shown later, Theorem \ref{theorem: extension of fundamental lemma and equivalence to SSKF} bears further implications including a handy closed-loop implementation of data-driven OP and the boundedness of the second moment of errors, which were not addressed by previous works.
\end{remark}

\begin{remark}  
    When innovations are not used, our data-driven OP reduces to the subspace predictive control (SPC) model \cite{favoreel1999spc}, where $g(t)$ in \eqref{equation: DeePC innovation division yf(t)} is given by:
    \begin{equation}
        \label{equation: SPC predictor}
        g_{\rm SPC}(t)=\mathbf{col}(U_{\rm p},U_{\rm f},Y_{\rm p})^\dagger \cdot \mathbf{col}(u_{\rm p}(t),u_{\rm f}(t),y_{\rm p}(t)).
    \end{equation}
    Even though the SPC model can feature asymptotic unbiasedness in the face of stochastic noise \cite{huang2008dynamic}, this critically relies on choosing a sufficiently large $L_{\rm p}$ such that $\bar{A}^{L_{\rm p}}\approx0$ and the bias can be eliminated. This, however, inevitably inflates the error variance \cite{breschi2023data,breschi2023tuning}. Conversely, choosing $L_{\rm p}$ small leads to a large bias, especially in a finite-sample regime. Nevertheless, this is not an issue for the proposed data-driven OP as it only requires $L_{\rm p}\ge\ell(A,B,C,D)$ in theory, which could be a key factor to its outperformance as to be shown in case studies.
\end{remark}

To further simplify the proposed data-driven OP, we may eliminate the equality constraint $\hat{E}_{\rm f}g(t)=0$ in \eqref{equation: DeePC innovation division g(t) with innovation estimates} by introducing $\hat{E}_{\rm f}^\bot$, i.e., the orthogonal projection matrix onto the null space of $\hat{E}_{\rm f}$. We assume that $\hat{E}_{\rm f}$ has full row rank, so the dimension of null space of $\hat{E}_{\rm f}$ is $N-L+1-L_{\rm f}n_y$. Thus, for each vector $g(t)\in\mathbb{R}^{N-L+1}$ in \eqref{equation: DeePC innovation division g(t) with innovation estimates}, there always exists a reduced-dimensional vector $h(t)\in\mathbb{R}^{N-L+1-L_{\rm f}n_y}$ satisfying $g(t)=\hat{E}_{\rm f}^\bot h(t),~ \hat{E}_{\rm f}\hat{E}_{\rm f}^\bot h(t)=0$. Then the predictor \eqref{equation: DeePC innovation division yf(t)} becomes $\hat{y}_{\rm f}(t)=Y_{\rm f}\hat{E}_{\rm f}^\bot h(t)$, where $h(t)$ satisfies the linear equations
\begin{equation}
    \label{equation: calculate alpha with precise innovation data}
    \begin{aligned}
        &\mathbf{col}(U_{\rm p}\hat{E}_{\rm f}^\bot,U_{\rm f}\hat{E}_{\rm f}^\bot,Y_{\rm p}\hat{E}_{\rm f}^\bot,\hat{E}_{\rm p}\hat{E}_{\rm f}^\bot)h(t)\\
        &=\mathbf{col}(u_{\rm p}(t),u_{\rm f}(t),y_{\rm p}(t),\hat{e}_{\rm p}(t)).   
    \end{aligned}
\end{equation}

In a similar spirit to \eqref{equation: pinv solution with estimated innovation}, a minimum-norm solution to \eqref{equation: calculate alpha with precise innovation data} is given by the Moore-Penrose inverse:
\begin{equation}
    \label{equation: minimum-norm solution alpha with precise innovation}
    \begin{split}
        h_{\rm pinv}(t)=&~\mathbf{col}(U_{\rm p}\hat{E}_{\rm f}^\bot,U_{\rm f}\hat{E}_{\rm f}^\bot,Y_{\rm p}\hat{E}_{\rm f}^\bot,\hat{E}_{\rm p}\hat{E}_{\rm f}^\bot)^\dagger \\
    &\quad \cdot \mathbf{col}(u_{\rm p}(t),u_{\rm f}(t),y_{\rm p}(t),\hat{e}_{\rm p}(t)).
    \end{split}
\end{equation}

\subsection{Closed-loop properties of data-driven OP}
We then delve into the closed-loop implementation of the proposed data-driven OP. Different from $u(t)$ and $y(t)$, the innovation $e(t)$ cannot be directly measured online. Thus, in order to implement the predictor in \eqref{equation: DeePC innovation division yf(t)} and \eqref{equation: DeePC innovation division g(t)} online, a key problem is how to derive $e(t)$ and update the vector of past innovations $e_{\rm p}(t)$. Interestingly, under the assumptions of Theorem \ref{theorem: extension of fundamental lemma and equivalence to SSKF}, the innovation at time $t$ can be readily derived from the data-driven OP as
\begin{equation}
\label{eq: computing innovation}
e(t)=y(t)-\hat{y}(t)=y(t)-\hat{y}_{{\rm f},1}(t),
\end{equation}
thanks to its equivalence with the SSKF-based OP. Hence, while moving on to the next instance $t+1$, the required vector $e_{\rm p}(t+1)$ can be updated in a moving window fashion:
\begin{equation}
    \label{equation: update of past innovation}
    \begin{split}
        e_{\rm p}(t+1) & = \mathbf{col}(e(t-L_{\rm p}+1),\cdots,e(t-1),e(t)) \\
        & = \mathbf{col}(e_{{\rm p},[2:L_{\rm p}]}(t),y(t)-\hat{y}_{{\rm f},1}(t)).
    \end{split}    
\end{equation}
In this sense, under the assumptions of Theorem \ref{theorem: extension of fundamental lemma and equivalence to SSKF}, the equivalence between the proposed data-driven OP and the SSKF-based OP remains valid in closed-loop. This inspires and lays the groundwork for using \eqref{eq: computing innovation} in the proposed data-driven OP based on innovation estimates. Another issue is how to initialize $e_{\rm p}(0)$ appropriately at the beginning, given $u_{\rm p}(0)$ and $y_{\rm p}(0)$. An intuitive idea is to leverage Theorem \ref{theorem: extension of fundamental lemma and equivalence to SSKF}(a) to attain $e_{\rm p}(0)$, whose rationale can be justified as follows.
\begin{clr}
\label{theorem: initialization of RHC scheme}
Under the assumptions behind Theorem \ref{theorem: extension of fundamental lemma and equivalence to SSKF}, $\{u_{\rm p}(t),y_{\rm p}(t),e_{\rm p}(t)\}$ is a valid trajectory of \eqref{equation: innovation-form LTI system} if and only if there exists a vector $g\in\mathbb{R}^{N-L+1}$ such that $\mathbf{col}(U_{\rm p},Y_{\rm p},E_{\rm p})g=\mathbf{col}(u_{\rm p}(t),y_{\rm p}(t),e_{\rm p}(t))$.
\end{clr}
\begin{proof}
    Since $\Tilde{u}^{\rm d}_{[1:N]}$ is persistently exciting of order $L+n_x$, $\mathcal{H}_{L+n_x}(\tilde{u}^{\rm d}_{[1:N]})$ has full row rank, so does the matrix $\mathcal{H}_{L_{\rm p}+n_x}(\tilde{u}^{\rm d}_{[1:N-L_{\rm f}]})$ constituted by rows of $\mathcal{H}_{L+n_x}(\tilde{u}^{\rm d}_{[1:N]})$. Thus, the truncated sequence $\tilde{u}^{\rm d}_{[1:N-L_{\rm f}]}$ is persistently exciting of order $L_{\rm p}+n_x$. Applying Theorem \ref{theorem: extension of fundamental lemma and equivalence to SSKF}(a) then completes the proof.
\end{proof}

Inspired by Corollary \ref{theorem: initialization of RHC scheme}, we propose to initialize $e_{\rm p}(0)$ through solving the following optimization problem:
\begin{equation}
    \label{equation: optimization problem for initial innovation data}
    \begin{aligned}
        \min_{e_{\rm p}(0),g}&\|e_{\rm p}(0)\|_2^2\\
        \rm{s.t.}~&\mathbf{col}(U_{\rm p},Y_{\rm p},E_{\rm p})g=\mathbf{col}(u_{\rm p}(0),y_{\rm p}(0),e_{\rm p}(0)).
    \end{aligned}
\end{equation}
which can be interpreted as a model-based moving horizon estimation (MHE) problem with known $(A, B, C, D, K)$:
\begin{equation}
    \label{equation: optimization problem for initial innovation data model-based}
    \begin{aligned}
        \min_{\hat{x}_{[-L_{\rm p}:0]},e}&~\|e\|_2^2\\
        \rm{s.t.}~~&~\mbox{Eq.}~\eqref{equation: innovation-form LTI system}, ~t=-L_{\rm p},...,-1,~e=e_{[-L_{\rm p}:-1]},
    \end{aligned}
\end{equation}
where the sum of one-step prediction errors is minimized. 
\begin{remark}
In \cite{wolff2024robust}, a data-driven MHE scheme was proposed for state estimation, where full state measurements are needed to construct Hankel matrices. Differently, the data-driven counterpart \eqref{equation: optimization problem for initial innovation data} of the model-based MHE \eqref{equation: optimization problem for initial innovation data model-based} requires no offline state measurements.
\end{remark}

Next, we present a formal analysis of closed-loop properties of data-driven OP $\hat{y}_{\rm f}(t)=Y_{\rm f}\hat{E}^\bot_{\rm f}h_{\rm pinv}(t)$ with $h_{\rm pinv}(t)$ given by \eqref{equation: minimum-norm solution alpha with precise innovation} that is built upon innovation estimates $\hat{e}^{\rm d}_{[1:N]}$ and a moving window update of its past innovations $\hat{e}_{\rm p}(t)$ in \eqref{equation: update of past innovation}. Our interest is concentrated on the boundedness of the second moment of its one-step ahead prediction error in closed-loop. To this end, we introduce an alternative ``implicit" data-driven OP, which is built upon ``true" values of innovation $e^{\rm d}_{[1:N]}$ and is known to coincide with the SSKF-based OP as per Theorem \ref{theorem: extension of fundamental lemma and equivalence to SSKF}(b). Thus, this ``implicit" OP always exists theoretically, which does not have to be physically implemented but is useful for further analysis. For clarity, we use the superscript $\cdot^{\rm KF}$ to indicate data matrices and variables pertaining to it. To describe the dynamics of $\hat{e}(t)$, we first delve into the gap between $\hat{y}(t)$ and $\hat{y}^{\rm KF}(t)\triangleq\hat{y}_{{\rm f},1}^{\rm KF}(t)$. Define $\alpha(t)=\hat{E}_{\rm f}^\bot h_{\mathrm{pinv}}(t)$ and $\alpha^{\rm KF}(t)=E_{\rm f}^\bot h^{\rm KF}_{\mathrm{pinv}}(t)$, where $h_{\rm pinv}^{\rm KF}(t)$ is the pseudo-inverse solution in the form of \eqref{equation: minimum-norm solution alpha with precise innovation} with $\{\hat{E}_{\rm p},\hat{E}_{\rm f}\}$ replaced by $\{E_{\rm p},E_{\rm f}\}$. It then follows that $\hat{y}^{\rm KF}(t) - \hat{y}(t)= Y_{{\rm f},[1:n_y]}[\alpha^{\rm KF}(t) - \alpha(t)]$. The dynamics of $\beta(t) = \alpha^{\rm KF}(t) - \alpha(t)$ is then described as follows.
\begin{theorem}
    \label{theorem: dynamic of Delta alpha}
    Under the assumptions in Theorem \ref{theorem: extension of fundamental lemma and equivalence to SSKF} and given an estimated sequence $\hat{e}^{\rm d}_{[1:N]}$, the dynamics of $\beta(t)$ in closed-loop is governed by:
    \begin{equation}
        \label{equation: dynamic of delta alpha}
        \begin{aligned}
            \beta(t+1)=&~MP\beta(t)+(M^{\rm KF}P^{\rm KF}-MP+\tilde{M} Z_1)\alpha^{\rm KF}(t)\\
            &+\tilde{M} Z_2u_{\rm f}(t+1)+\tilde{M} Z_3e(t)
        \end{aligned}
    \end{equation}
    where
    \begin{equation*}
        \begin{aligned}
        &M^{\rm KF}=E_{\rm f}^\bot\mathbf{col}(U_{\rm p}E_{\rm f}^\bot,U_{\rm f}E_{\rm f}^\bot,Y_{\rm p}E_{\rm f}^\bot,E_{\rm p}E_{\rm f}^\bot)^\dagger,\\
        &P^{\rm KF}=\mathbf{col}(U_{{\rm p},[n_u+1:n_uL_{\rm p}]},U_{{\rm f},[1:n_u]},0_{n_uL_{\rm f}\times(N-L+1)},\\
        &\quad Y_{{\rm p},[n_y+1:n_yL_{\rm p}]},0_{n_y\times(N-L+1)},E_{{\rm p},[n_y+1:n_yL_{\rm p}]},-Y_{{\rm f},[1:n_y]}),\\
        &Z_1=\mathbf{col}(0_{[n_uL+n_y(L_{\rm p}-1)]\times(N-L+1)},Y_{{\rm f},[1:n_y]},\\
        &\quad\quad\quad\quad\quad\quad\quad0_{n_y(L_{\rm p}-1)\times(N-L+1)},Y_{{\rm f},[1:n_y]}),\\
        &Z_2=\mathbf{col}(0_{n_uL_{\rm p}\times n_uL_{\rm f}},I_{n_uL_{\rm f}},0_{2n_yL_{\rm p}\times n_uL_{\rm f}}),\\
        &Z_3=\mathbf{col}(0_{[n_uL+n_y(L_{\rm p}-1)]\times n_y},I_{n_y},0_{[n_y(L_{\rm p}-1)]\times p},I_{n_y}),
        \end{aligned}
    \end{equation*}
    matrices $M$ and $P$ are defined akin to $M^{\rm KF}$ and $P^{\rm KF}$ by replacing $\{E_{\rm p},E_{\rm f}\}$ with $\{\hat{E}_{\rm p},\hat{E}_{\rm f}\}$, and $\tilde{M}=M^{\rm KF}-M$.
\end{theorem}
\begin{proof}
    We rewrite $\alpha^{\rm KF}(t+1)$ as:
    \begin{equation}
        \label{equation: minimum-norm solution at t+1 with KF}
        \begin{aligned}
            &\alpha^{\rm KF}(t+1)\\
            =&M^{\rm KF}\mathbf{col}(u_{\rm p}(t+1),u_{\rm f}(t+1),y_{\rm p}(t+1),e_{\rm p}(t+1)).
        \end{aligned}
    \end{equation}
    For the ``implicit'' data-driven OP equal to the SSKF-based predictor, it holds that:
    \begin{equation}
        \label{equation: calculate alpha with precise innovation data KF}
        \begin{aligned}
            &\mathbf{col}(U_{\rm p}E_{\rm f}^\bot,U_{\rm f}E_{\rm f}^\bot,Y_{\rm p}E_{\rm f}^\bot,E_{\rm p}E_{\rm f}^\bot)h(t)\\
            &=\mathbf{col}(u_{\rm p}(t),u_{\rm f}(t),y_{\rm p}(t),e_{\rm p}(t)),   
        \end{aligned}
    \end{equation}
    which is similar to \eqref{equation: calculate alpha with precise innovation data}. Based on \eqref{equation: calculate alpha with precise innovation data KF}, one obtains:    
    \begin{equation*}
        \begin{aligned}
            &u_{\rm p}(t+1)=\begin{bmatrix}u_{{\rm p},[2:L_{\rm p}]}(t)\\u(t)\end{bmatrix}
            =\begin{bmatrix}
                U_{{\rm p},[n_u+1:n_uL_{\rm p}]}\\U_{{\rm f},[1:n_u]}
            \end{bmatrix}\alpha^{\rm KF}(t),\\
            &y_{\rm p}(t+1)=\begin{bmatrix}y_{{\rm p},[2:L_{\rm p}]}(t)\\\hat{y}^{\rm KF}(t)+e(t)\end{bmatrix}
            =\begin{bmatrix}
                Y_{{\rm p},[n_y+1:n_yL_{\rm p}]}\alpha^{\rm KF}(t)\\
                Y_{{\rm f},[1:n_y]}\alpha^{\rm KF}(t)+e(t)
            \end{bmatrix},\\
            &e_{\rm p}(t+1)=\begin{bmatrix}e_{{\rm p},[2:L_{\rm p}]}(t)\\e(t)\end{bmatrix}
            =\begin{bmatrix}
                E_{{\rm p},[n_y+1:n_yL_{\rm p}]}\alpha^{\rm KF}(t)\\
                e(t)
            \end{bmatrix},
        \end{aligned}
    \end{equation*}
    where we utilize the fact that $y(t) = \hat{y}^{\rm KF}(t) + e(t)$ and $\hat{y}^{\rm KF}(t) = Y_{{\rm f},[1:n_y]}\alpha^{\rm KF}(t)$. Plugging them into \eqref{equation: minimum-norm solution at t+1 with KF} yields:
    \begin{equation}
        \label{equation: dynamic of alpha}
        \begin{aligned}
            &\alpha^{\rm KF}(t+1)\\
            =&~ M^{\rm KF}[(P^{\rm KF}+Z_1)\alpha^{\rm KF}(t)+Z_2u_{\rm f}(t+1)+Z_3e(t)].
        \end{aligned}
    \end{equation}
    Meanwhile, $\alpha(t+1)$ can be rewritten as:
    \begin{equation}
        \label{equation: minimum-norm solution at t+1}
        \begin{aligned}
            &\alpha(t+1)\\
            =&M\mathbf{col}(u_{\rm p}(t+1),u_{\rm f}(t+1),y_{\rm p}(t+1),\hat{e}_{\rm p}(t+1)).
        \end{aligned}
    \end{equation}
    Using \eqref{equation: calculate alpha with precise innovation data}, one obtains:
    \begin{equation*}
        \begin{aligned}
            &u_{\rm p}(t+1)=\mathbf{col}(U_{{\rm p},[n_u+1:n_uL_{\rm p}]},U_{{\rm f},[1:n_u]})\alpha(t),\\
            &y_{\rm p}(t+1)=\mathbf{col}(Y_{{\rm p},[n_y+1:n_yL_{\rm p}]}\alpha(t),Y_{{\rm f},[1:n_y]}\alpha^{\rm KF}(t)+e(t)),\\
            &\hat{e}_{\rm p}(t+1)=\mathbf{col}(\hat{E}_{{\rm p},[n_y+1:n_yL_{\rm p}]}\alpha(t),Y_{{\rm f},[1:n_y]}\beta(t)+e(t)),
        \end{aligned}
    \end{equation*}
    where the last element of $\hat{e}_{\rm p}(t+1)$, namely $\hat{e}(t)$, is expressed as:
    \begin{equation}
        \label{equation: one-step prediction error with imprecise innovation}
        \begin{aligned}
        \hat{e}(t)=y(t)-\hat{y}(t)
        &=y(t)-\hat{y}^{\rm KF}(t)+\hat{y}^{\rm KF}(t)-\hat{y}(t)\\
        &=e(t)+Y_{{\rm f},[1:n_y]}\beta(t).
        \end{aligned}
    \end{equation}
    By plugging $\{u_{\rm p}(t+1),y_{\rm p}(t+1),\hat{e}_{\rm p}(t+1)\}$ into \eqref{equation: minimum-norm solution at t+1}, one obtains:
    \begin{equation}
        \label{equation: dynamic of alpha1}
        \begin{split}
        &\alpha(t+1) \\
        =&M[P\alpha(t)+Z_1\alpha^{\rm KF}(t)+Z_2u_{\rm f}(t+1)+Z_3e(t)].
        \end{split}
    \end{equation}
    Combining \eqref{equation: dynamic of alpha} and \eqref{equation: dynamic of alpha1} yields \eqref{equation: dynamic of delta alpha}.
\end{proof}
Now we arrive at the following result on the second moment boundedness of prediction errors.
\begin{theorem}
    \label{theorem: second moment bounded}
    Let the assumptions in Theorem \ref{theorem: extension of fundamental lemma and equivalence to SSKF} hold and an innovation estimate sequence $\hat{e}^{\rm d}_{[1:N]}$ be given. If the control input $u(t)$ is bounded and ensures the boundedness of $\mathbb{E}[\|y(t)\|^2]$, then the following statements hold.
    \begin{enumerate}
    \item[(a)] The second moment of $\alpha^{\rm KF}(t)$ is bounded.
    \item[(b)] If $\Theta=MP$ is Schur stable, then $\beta(t)$ has a bounded second moment, and so is the one-step prediction error $\hat{e}(t)$ of the proposed data-driven OP based on innovation estimates $\hat{e}^{\rm d}_{[1:N]}$.
    \end{enumerate}
\end{theorem}
\begin{proof}
    Since $e(t)$ is the one-step prediction error from SSKF \eqref{equation: steady-state kalman filter}, $\mathbb{E}[\|e_{\rm p}(t)\|^2]$ is bounded \cite{huang2008dynamic}. Because $u(t)$ and $\mathbb{E}[\|y(t)\|^2]$ are also bounded, the boundedness of $\mathbb{E}[\|h^{\rm KF}_{\rm{pinv}}(t)\|^2]$ then follows from \eqref{equation: minimum-norm solution alpha with precise innovation}, and thus the boundedness of $\mathbb{E}[\|\alpha^{\rm KF}(t)\|^2]$ can be proven by definition. Next, we cope with (b). Note that $\beta(t)$ is the system state of \eqref{equation: dynamic of delta alpha}, where the inputs consist of $\alpha^{\rm KF}(t)$, $u_{\rm f}(t+1)$ and $e(t)$, all of which have bounded second moments. Thus, $\mathbb{E}[\|\beta(t)\|^2]$ is bounded due to the Schur stability of $\Theta$. Then, according to \eqref{equation: one-step prediction error with imprecise innovation}, the one-step prediction error $\hat{e}(t)$ induced by the proposed OP has a bounded second moment due to the boundedness of $\mathbb{E}[\|e(t)\|^2]$ and $\mathbb{E}[\|\beta(t)\|^2]$.
\end{proof}
\begin{remark}
    In Theorem \ref{theorem: second moment bounded}, the expectation at time $t$ is taken with respect to the joint distribution of $x(0)$, $w_{[0:t-1]}$ and $v_{[0:t]}$ in \eqref{equation: LTI system}, which are also the source of uncertainty for the innovation form \eqref{equation: innovation-form LTI system}. To ensure boundedness of $\mathbb{E}[\|y(t)\|^2]$, $u(t)$ can be either stochastic, e.g. from a feedback control law, or deterministic. In the former case, the randomness of $u(t)$ stems from $x(0)$, $w_{[0:t-1]}$ and $v_{[0:t-1]}$, so the proof of Theorem \ref{theorem: second moment bounded} still holds. In the latter case, Theorem \ref{theorem: second moment bounded} implicitly requires \eqref{equation: LTI system} to be open-loop stable and does not apply to unstable open-loop systems with deterministic inputs.
\end{remark}
\begin{remark}
    Indeed, the conditions in Theorem \ref{theorem: second moment bounded} are also sufficient for establishing the boundedness of the second moment of multi-step prediction error $\tilde{y}_{\rm f}(t)=y_{\rm f}(t)-\hat{y}_{\rm f}(t)$, where $y_{\rm f}(t)=y_{[t:t+L_{\rm f}-1]}$. Similar to \eqref{equation: one-step prediction error with imprecise innovation}, $\tilde{y}_{\rm f}(t)$ can be decomposed as $\tilde{y}_{\rm f}(t)=y_{\rm f}(t)-\hat{y}_{\rm f}(t)=y_{\rm f}(t)-\hat{y}^{\rm KF}_{\rm f}(t)+Y_{\rm f}\beta(t)$. Then the boundedness of $\mathbb{E}[\|\tilde{y}_{\rm f}(t)\|^2]$ can be established akin to Theorem \ref{theorem: second moment bounded}.
\end{remark}
Note that the construction of $\Theta$ only relies on the pre-collected trajectory $\{u^{\rm d}(t),y^{\rm d}(t)\}_{i=1}^N$ and the innovation estimates $\hat{e}^{\rm d}_{[1:N]}$. The proof of Theorem \ref{theorem: second moment bounded} indicates that when $\Theta$ bears no Schur stability, the prediction error $\hat{e}(t)$ can diverge, which invalidates the data-driven OP. Thus, the stability of $\Theta$ is a practical indicator for the ``validity'' of innovation estimates $\hat{e}^{\rm d}_{[1:N]}$. As such, the stability of $\Theta$ shall be examined as a key step prior to the online adoption of the data-driven OP. The implementation procedure is detailed in Algorithm \ref{algorithm: receding horizon simulation scheme}. Specifically, whenever instability of $\Theta$ is seen after deriving $\hat{e}^{\rm d}_{[1:N]}$, one may adjust the choice of $\rho$ in VARX modelling, or collect another trajectory $\{u^{\rm d}(t),y^{\rm d}(t)\}_{i=1}^N$ and estimate innovations to attain a stable $\Theta$.

\begin{algorithm}[ht]
    \caption{Implementation of Data-Driven OP}
    \label{algorithm: receding horizon simulation scheme}
    \textbf{Offline data collection and pre-processing:}
    \begin{algorithmic}[1]
        \State Collect data $\{u^{\rm d}(i),y^{\rm d}(i)\}_{i=-N_\rho}^N$;
        \State Choose $\rho$ such that $1\le\rho\le N_\rho$;
        \State Estimate $\hat{e}^{\rm d}_{[1:N]}$ by VARX modelling;
        \State If $\Theta$ is Schur stable, go to Step 5, else choose another $1\le\rho\le N_{\rho}$ and go to Step 3, or directly go to Step 1;
        \State Construct $\{U_{\rm d},Y_{\rm d},\hat{E}_{\rm d}\}$ from $\{u^{\rm d}(i),y^{\rm d}(i),\hat{e}^{\rm d}(i)\}_{i=1}^N$;
        \State Collect an initial past trajectory $\{u(i),y(i)\}_{i=-L_{\rm p}}^{-1}$ and initialize $\hat{e}_{\rm p}(0)$ based on \eqref{equation: optimization problem for initial innovation data};
    \end{algorithmic}
    \textbf{Online output prediction:}
    \begin{algorithmic}[1]
        \State Given a future input trajectory $u_{\rm f}(t)=u_{[t:t+L_{\rm f}-1]}$, derive future output prediction $\hat{y}_{\rm f}(t)=Y_{\rm f}\hat{E}_{\rm f}h_{\rm pinv}(t)$ with $h_{\rm pinv}(t)$ given by \eqref{equation: minimum-norm solution alpha with precise innovation};
        \State Collect $y(t)$;
        \State Update $\hat{e}_{\rm p}(t+1)$ based on \eqref{equation: update of past innovation} and update $\{u_{\rm p}(t+1), y_{\rm p}(t+1)\}$;
        \State $t\gets t+1$, go to Step 1.
    \end{algorithmic}
\end{algorithm}

\begin{remark}
    When $\Theta$ is Schur stable, the second moment of prediction error $\hat{e}(t)$ is known to be bounded. As such, we can run the data-driven OP on a new validation dataset to collect empirical samples of prediction errors, based on which statistics of $\hat{e}(t)$ such as the covariance can be computed to characterize the uncertainty.
\end{remark}

\subsection{Application to data-driven predictive control}

In the following, we discuss the application of the proposed data-driven OP to predictive control tasks. The trivial inclusion of data-driven OP $Y_{\rm f}g(t)=y_{\rm f}(t)$ and \eqref{equation: DeePC division} gives rise to the well-known formulation of DeePC \cite{coulson2019data}:
\begin{equation}
    \label{eq: DeePC}
    \begin{aligned}
        \min_{u_{\rm f}(t),\hat{y}_{\rm f}(t),g(t)}~~&\mathcal{J}(u_{\rm f}(t),\hat{y}_{\rm f}(t))\\
        {\rm s.t.}~\quad\quad&\hat{y}_{\rm f}(t)=Y_{\rm f}g(t),~\mbox{Eq.}~\eqref{equation: DeePC division},\\
        &u_{\rm f}(t)\in\mathbb{U},~\hat{y}_{\rm f}(t)\in\mathbb{Y},
    \end{aligned}
\end{equation}
where $\mathbb{U} = \mathbb{U}_1 \times \cdots \times \mathbb{U}_{L_{\rm f}}$, $\mathbb{Y} = \mathbb{Y}_1 \times \cdots \times \mathbb{Y}_{L_{\rm f}}$ are Cartesian products of sets with stage-wise constraints $\mathbb{U}_k$ and $\mathbb{Y}_k$. Note that $g(t)$ appears as a decision variable that shall be optimized together with $u_{\rm f}(t)$ and $\hat{y}_{\rm f}(t)$. The objective function of \eqref{eq: DeePC} can be chosen as the standard quadratic cost $\mathcal{J}(u_{\rm f}(t),\hat{y}_{\rm f}(t)) = \sum_{\tau=t}^{t+L_{\rm f}-1} ||\hat{y}(\tau)-r(\tau)||_{Q}^2+||{u}(\tau)||_{R}^2$, where $Q,R \succ 0$ are cost weighting matrices and $r(t)$ is a reference signal. Alternatively, $\mathcal{J}(\cdot, \cdot)$ may also be designed to track a desired equilibrium \cite[Eq. (2a)]{berberich2021data} or a reference input-output sequence \cite[Eq. (2)]{dorfler2022bridging}. To hedge against inflated variance of OP under uncertainty, a common option is to use the pseudo-inverse solution \eqref{equation: SPC predictor}, yielding the generic formulation of SPC \cite{favoreel1999spc}:
\begin{equation}
    \label{SPC opt problem}
    \begin{aligned}
    \min_{u_{\rm f}(t),\hat{y}_{\rm f}(t)}\quad&\mathcal{J}(u_{\rm f}(t),\hat{y}_{\rm f}(t)) \\
    \mathrm{s.t.}~~~\quad&\hat{y}_{\rm f}(t)=Y_{\rm f}g_{\rm SPC}(t),~\mbox{Eq.}~\eqref{equation: SPC predictor},\\
    &u_{\rm f}(t)\in\mathbb{U},~\hat{y}_{\rm f}(t)\in\mathbb{Y}.
    \end{aligned}
\end{equation}
By the same token, one can also insert the proposed data-driven OP $\hat{y}_{\rm f}(t)=Y_{\rm f}\hat{E}_{\rm f}h_{\rm pinv}(t)$ with $h_{\rm pinv}(t)$ given by \eqref{equation: minimum-norm solution alpha with precise innovation} into predictive control problems. In this way, we arrive at Inno-DeePC, a new innovation-based data-driven predictive control formulation:
\begin{equation}
    \label{equation: SPC-IV}
    \begin{aligned}
        \min_{u_{\rm f}(t),\hat{y}_{\rm f}(t)}\quad&\mathcal{J}(u_{\rm f}(t),\hat{y}_{\rm f}(t))\\
    \rm{s.t.}~\quad~~&\hat{y}_{\rm f}(t)=Y_{\rm f}\hat{E}_{\rm f}h_{\rm pinv}(t),~\mbox{Eq.}~\eqref{equation: minimum-norm solution alpha with precise innovation},\\
    &u_{\rm f}(t)\in\mathbb{U},~\hat{y}_{\rm f}(t)\in\mathbb{Y}.
    \end{aligned}
\end{equation}
In practice, the data-driven OP used in \eqref{equation: SPC-IV} builds upon innovation estimates $\hat{e}^{\rm d}_{[1:N]}$ and thus bypasses the usage of a parametric model. Recall that under the assumptions in Theorem \ref{theorem: extension of fundamental lemma and equivalence to SSKF}, e.g. when the true trajectory of innovations is available, the data-driven OP in \eqref{equation: SPC-IV} becomes equal to the SSKF-based OP. In this sense, the enabled control design \eqref{equation: SPC-IV} can be interpreted as a data-driven realization of the model-based predictive control scheme where an SSKF-based OP \eqref{equation: model-based simulation} is adopted for multi-step prediction.

\begin{remark}
    The proposed Inno-DeePC is reminiscent of the known DeePC with Extended Kalman Filter (EKF) in \cite{alpago2020extended}, where an EKF is adopted to hedge against noise in past outputs $y_{\rm p}(t)$. Notably, in \cite{alpago2020extended} the state variables of EKF are chosen as $y_p(t)$ rather than state estimates of the LTI system \eqref{equation: LTI system}. Differently, the proposed Inno-DeePC in \eqref{equation: SPC-IV} can be viewed as implicitly estimating states of \eqref{equation: LTI system} for output prediction, owing to a data-driven (albeit approximated) implementation of a model-based SSKF.
\end{remark}

\section{Numerical Examples}
In this section, numerical simulations are carried out to empirically investigate output prediction and control performance of the proposed methods based on innovation estimates. Consider the LTI system \eqref{equation: LTI system} with $A=\begin{bmatrix}0.7326&-0.0861\\0.1722&0.9909\end{bmatrix}
        ,~
        B = \begin{bmatrix}
            0.0609\\0.0064
        \end{bmatrix},~
        C=\begin{bmatrix}
            0&1.4142
        \end{bmatrix}
        ,~
        D = 0$ \cite{breschi2023data},
where $\{w(t),v(t)\}$ are set to be zero-mean Gaussian distributed with $\Sigma_w=\sigma_w^2I_2$ and $\Sigma_v=\sigma_v^2$. Here, we set $\sigma_w=q\times 10^{-4}$ and $\sigma_v=4.5q\times10^{-4}$ such that the noise level can be adjusted via a single parameter $q>0$.

\subsection{Results of data-driven output predictor}

We first investigate the performance of the following three predictors.
\begin{enumerate}
    \item \textbf{Inno-OP}: The proposed data-driven OP $\hat{y}_{\rm f}(t)=Y_{\rm f}\hat{E}_{\rm f}h_{\rm pinv}(t)$ with $h_{\rm pinv}(t)$ given by \eqref{equation: minimum-norm solution alpha with precise innovation}, as implemented in Algorithm \ref{algorithm: receding horizon simulation scheme}.
    \item \textbf{PBSID}: The model-based OP \eqref{equation: steady-state kalman filter} using system matrices $(A,B,C,D,K)$ identified by predictor-based subspace identification (PBSID) \cite{houtzager2009varmax}.
    \item \textbf{SPC}: The SPC model \cite{favoreel1999spc}, i.e. $\hat{y}_{\rm f}(t)=Y_{\rm f}g_{\rm SPC}(t)$ with $g_{\rm SPC}(t)$ given by \eqref{equation: SPC predictor}.
\end{enumerate}

\begin{figure}[ht]
    \centering
    \includegraphics[width=0.48\textwidth]{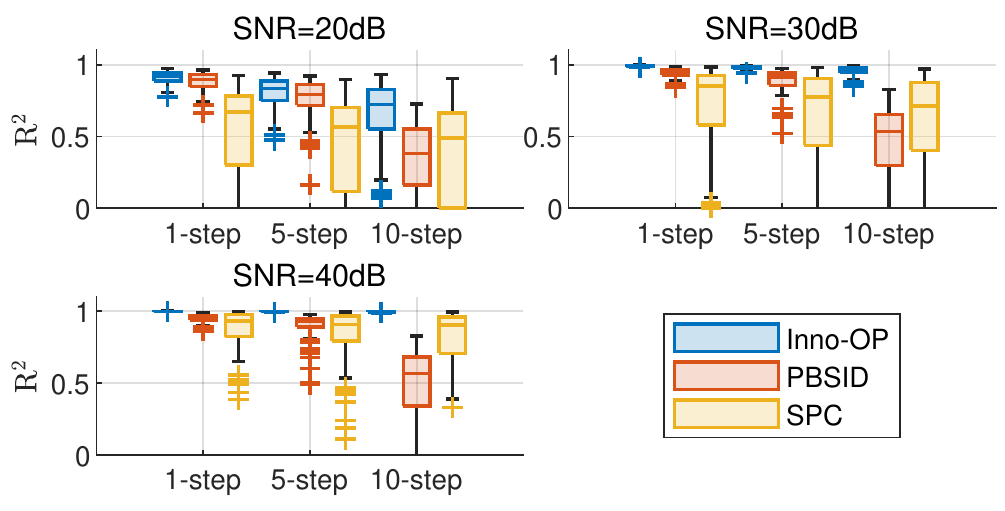}
     \caption{Boxplots of ${\rm R}^2$ indices of different multi-step prediction at different SNRs in $100$ simulations.}
    \label{fig:prediction performance}
\end{figure}

We set $L_{\rm p}=10$, which ensures $\|\bar{A}^{L_{\rm p}}\|_2\le0.005$ to reduce the bias in SPC (see Remark 4) and also satisfies the requirement for Inno-DeePC. Meanwhile, $L_{\rm f}=15$. A square wave with a period of $50$ time-steps and amplitude of $2$, contaminated by a zero-mean Gaussian distributed sequence with variance $0.01$, is used as the offline input $u^{\rm d}_{[1:N]}$. Based on this, an input-output trajectory of length $N=200$ is pre-collected to build various OPs. To evaluate their prediction performance, we generate a test trajectory of length $N_{\rm test}=100$, where the inputs are sampled from a zero-mean Gaussian distribution with variance $4$. For a comprehensive assessment of prediction performance, 100 Monte Carlo runs are carried out at different signal-to-noise ratios (SNRs)\footnote{We define ${\rm SNR}=10\log_{10} \frac{{\rm var}[y^{\rm d}(t)-e^{\rm d}(t)]}{{\rm var}[e^{\rm d}(t)]}$, where $e^{\rm d}(t)$ is derived by running the model-based SSKF and $K$ is obtained by solving the Riccati equation.}, i.e., SNR = $20$, $30$, $40$dB, by setting $q=11.49,~1.13,~0.11$, respectively. In each Monte Carlo run, both offline input-output trajectory and test trajectory are randomly generated at the same SNR. To evaluate the prediction accuracy, the coefficient of determination ${\rm R^2}=1- \sum_{k=1}^{N_{\rm test}}[y(k)-\hat{y}(k)]^2/\sum_{k=1}^{N_{\rm test}}[y(k)-\bar{y}]^2$ is adopted, where $\bar{y}$ denotes the average of $y(k)$. Fig. \ref{fig:prediction performance} depicts the ${\rm R^2}$ indices of 1-step, 5-step, and 10-step predictions of different predictors. It can be seen that the prediction performance of all methods degrades as SNR decreases and the prediction step becomes longer. It is obvious that the proposed Inno-OP obtains the highest prediction accuracy under all circumstances, and its advantage becomes more pronounced for longer predictions. Moreover, the values of ${\rm R}^2$ indicate that Inno-OP yields predictions that can well explain the total variability of outputs, thereby showing promise in the face of data uncertainty.

\begin{figure}[ht]
    \centering
    \includegraphics[width=0.35\textwidth]{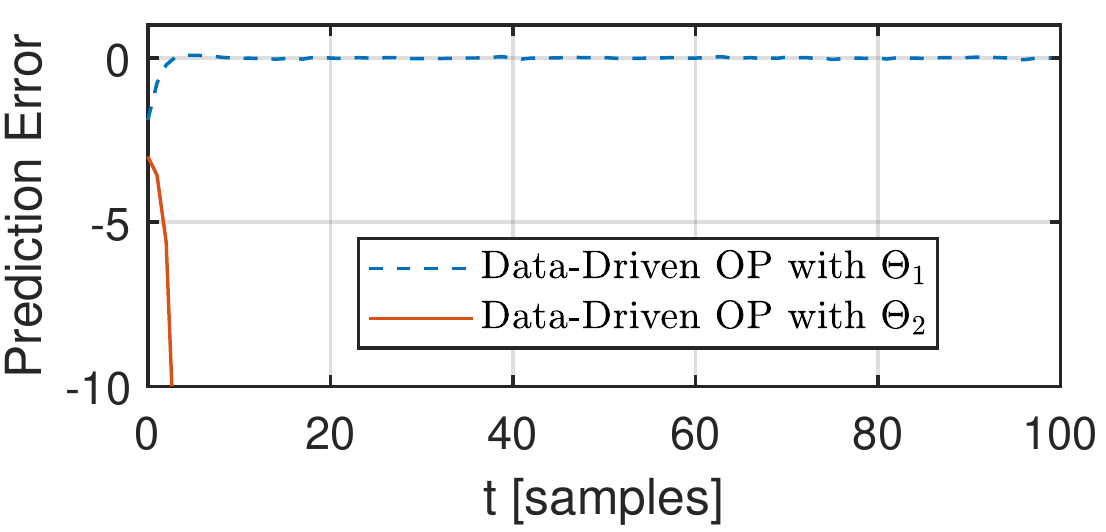}
    \caption{One-step prediction errors of data-driven OPs with $\Theta_1$ and $\Theta_2$ in a single Monte Carlo run.}
    \label{fig:CEOP unstable}
\end{figure}

Then we showcase the correctness of Theorem \ref{theorem: second moment bounded} based on a single Monte Carlo run with ${\rm SNR}=\SI{30}{dB}$. We create two trajectories of innovation estimates by varying $\rho$ and implement the resultant data-driven OP in Algorithm 1. Specifically, using $\rho=15$ yields a stable matrix $\Theta_1$, whereas using $\rho=50$ yields an unstable matrix $\Theta_2$.\footnote{The instability of $\Theta_2$ is mainly due to overfitting of VARX caused by setting a large $\rho=50$, as discussed in Appendix.} In this way, there are two data-driven OPs developed and their one-step prediction errors are profiled in Fig. \ref{fig:CEOP unstable}. Clearly, prediction errors induced by $\Theta_1$ show a bounded variance, while those induced by $\Theta_2$ diverge, thereby clearly showing the correctness of Theorem \ref{theorem: second moment bounded}.

\subsection{Results of data-driven predictive control}

Next, we investigate the closed-loop control performance of Inno-DeePC based on innovation estimates. For comparison, five control strategies are employed.
\begin{enumerate}
    \item \textbf{SSKF-MPC}: The model-based ``oracle" control scheme with the SSKF-based OP \eqref{equation: model-based simulation}, where $(A,B,C,D,K)$ is assumed to be known and the SSKF is used for prediction. 
    \item \textbf{Inno-DeePC}: The proposed Inno-DeePC scheme \eqref{equation: SPC-IV} based on the data-driven OP in Algorithm \ref{algorithm: receding horizon simulation scheme}, where $\hat{e}_{\rm p}(0)$ is initialized by solving \eqref{equation: optimization problem for initial innovation data} with $\hat{E}_{\rm p}$. 
    \item \textbf{PCE-DeePC}: The stochastic DeePC scheme based on polynomial chaos expansions \cite{pan2021stochastic}.
    \item \textbf{SPC}: The classical SPC scheme \eqref{SPC opt problem} \cite{favoreel1999spc}.
    \item \textbf{Reg-DeePC}: The regularized DeePC \cite[Eq. (16)]{dorfler2022bridging}, where the regularization parameter $\lambda$ is selected from a grid of values within $[10^{-2},10^4]$.
\end{enumerate}

\begin{figure}[ht]
    \centering
    \subfigure[SNR=$20$dB]{\includegraphics[width=0.23\textwidth]{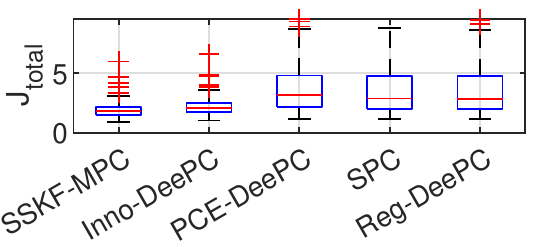}}
    \subfigure[SNR=$30$dB]{\includegraphics[width=0.23\textwidth]{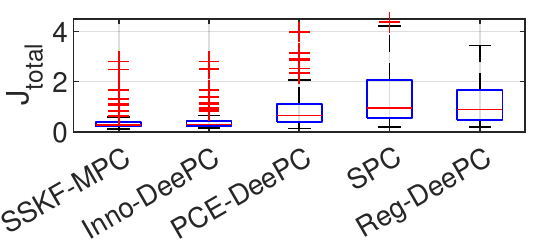}}
    \subfigure[SNR=$40$dB]{\includegraphics[width=0.23\textwidth]{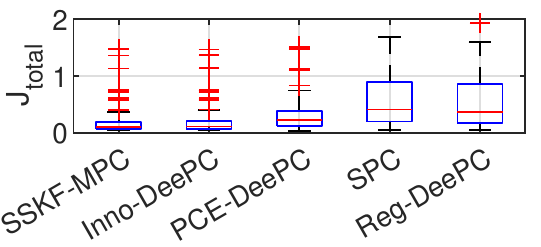}}
     \caption{Boxplots of $\mathcal{J}_{\rm total}$ at different SNRs in $100$ simulations.}
    \label{fig: control performance}
\end{figure}

The experimental settings are identical to those in the previous subsection. For predictive control design, we choose $Q=1$, $R=0.01$, the reference $r(k)=\mathrm{sin}(2\pi k/N_{\rm test}),~k=1,2,...,N_{\rm test}$ and the stage-wise constraints as $\mathbb{U}_k=\{u|-2\le u\le2\}$ and $\mathbb{Y}_k=\{\hat{y}|-2\le\hat{y}\le2\}$. 
To evaluate the control performance, the indices $\mathcal{J}_y = \sum_{k=1}^{N_{\rm test}}||y(k) - r(k)||_Q^2$, $\mathcal{J}_u = \sum_{k=1}^{N_{\rm test}}||u(k)||_R^2$ and $\mathcal{J}_{\rm total}=\mathcal{J}_y + \mathcal{J}_u$ are used, where $u(k)$ is the control action decided by each control strategy and $y(k)$ denotes the realized system output. The boxplots of $\mathcal{J}_{\rm total}$ in $100$ Monte Carlo simulations are shown in Fig. \ref{fig: control performance}, and more detailed results of $\mathcal{J}_u$ and $\mathcal{J}_y$ are presented in Table \ref{tab:predictive control performance}. It can be seen that PCE-DeePC achieves better performance than the generic SPC and Reg-DeePC. Meanwhile, the proposed Inno-DeePC outperforms PCE-DeePC, SPC, and Reg-DeePC under all settings, featuring a narrower gap with SSKF-MPC. 
This demonstrates that the proposed Inno-DeePC scheme draws strength from using innovations in improving the control performance of stochastic systems.

\begin{table*}
    \centering
    \caption{Statistics of $\mathcal{J}_u$ and $\mathcal{J}_y$ (Mean $\pm$ Standard Deviation) at Different SNRs in $100$ Monte Carlo Simulations}
    \label{tab:predictive control performance}
    \begin{tabular}{ccccccc}
    \toprule
    &\multicolumn{2}{c}{SNR = 20dB}&\multicolumn{2}{c}{SNR = 30dB}&\multicolumn{2}{c}{SNR = 40dB}\\
    \cmidrule{2-7}
    &$\mathcal{J}_u$&$\mathcal{J}_y$&$\mathcal{J}_u$&$\mathcal{J}_y$&$\mathcal{J}_u$&$\mathcal{J}_y$\\
    \midrule
    SSKF-MPC&$1.11\pm0.13$&$1.93\pm0.75$& $0.83\pm0.07$&$0.38\pm0.42$&$0.80\pm0.06$&$0.22\pm0.38$\\
    Inno-DeePC&$1.04\pm0.15$&$2.23\pm0.82$&$0.84\pm0.07$&$0.41\pm0.42$&$0.80\pm0.06$&$0.23\pm0.38$\\
    PCE-DeePC&$1.10\pm0.31$&$4.44\pm4.61$&$0.87\pm0.17$&$0.88\pm0.80$&$0.80\pm0.10$&$0.33\pm0.39$\\
    SPC&$1.10\pm0.37$&$5.26\pm9.21$&$0.94\pm0.32$&$2.40\pm5.03$&$0.84\pm0.18$&$0.68\pm0.81$\\
    Reg-DeePC&$1.02\pm0.33$&$4.06\pm3.70$&$0.84\pm0.20$&$1.67\pm2.87$&$0.83\pm0.17$&$0.65\pm0.79$\\
    \bottomrule
    \end{tabular}
\end{table*}

\section{Concluding Remarks}
We proposed a new data-driven approach to output prediction and control of stochastic LTI systems. Using the innovation form, all additive uncertainty is recaptured in terms of innovations that can be readily estimated from input-output data. By applying the fundamental lemma, we derived an innovation-based data-driven OP. Its equivalence with the SSKF-based OP was established, inspiring an easy closed-loop implementation. We also presented a sufficient condition for the boundedness of the second moment of prediction errors. Using the proposed data-driven OP in control design gives rise to a new DeePC formulation. Numerical case studies showed the performance improvement of the proposed methods over existing formulations. In future work, a direction is to consider stochastic control problems, e.g. chance-constrained control, to handle the uncertainty of OP.

\section*{Appendix: Data-Driven Innovation Estimation}
Due to the Schur stability of $\bar{A} = A-KC$, $\bar{A}^\rho\approx0$ holds for a sufficiently large $\rho \ge 0$. Thus, the VARX model can be compactly expressed in a matrix form \cite{chiuso2007role}:
\begin{equation*}
\begin{aligned}
    \mathcal{H}_1(y_{[1:N]})
    \approx &~ \Phi_y \mathcal{H}_\rho(y_{[1-\rho:N-1]}) + \Phi_u \mathcal{H}_\rho(u_{[1-\rho:N-1]}) \\
    &~ + D\mathcal{H}_1(u_{[1:N]}) + \mathcal{H}_1(e_{[1:N]}),
\end{aligned}  
\label{eq:41b}
\end{equation*}
where $\Phi_y \triangleq C[\bar{A}^{\rho-1}K~\cdots~\bar{A}K~K]$ and $\Phi_u \triangleq C[\bar{A}^{\rho-1}\bar{B}~\cdots~\bar{A}\bar{B}~\bar{B}]$ enclose MPs of \eqref{equation: innovation-form LTI system} with $\bar{B}\triangleq B-KD$. Solving the following least-squares regression (LSR) yields estimates of MPs $\{\Phi_y, \Phi_u \}$ and $D$ \cite{van2013closed}:
\begin{equation}
    \label{eq: LSR}
    \min_{\Phi_y, \Phi_u, D} \left \| \mathcal{H}_1(y_{[1:N]}) - 
    \begin{bmatrix}
        \Phi_y~ \Phi_u~ D
    \end{bmatrix}\Upsilon \right \|_F^2,
\end{equation}
where $\Upsilon=\mathbf{col}(\mathcal{H}_\rho(y_{[1-\rho:N-1]}),\mathcal{H}_\rho(u_{[1-\rho:N-1]}),\mathcal{H}_1(u_{[1:N]}))$ involves additional input-output data $\{u(i),y(i)\}_{i=1-\rho}^0$ preceding $\{u(i),y(i)\}_{i=1}^N$. Note that the residuals of \eqref{eq: LSR} provide an explicit approximation of $\mathcal{H}_1(e_{[1:N]})$, based on which $\hat{E}_d$ can be constructed. Besides, $\rho$ shall be suitably chosen to achieve a trade-off. Specifically, a large $\rho > n$ is needed to reduce approximation errors caused by $\bar{A}^\rho$, but choosing $\rho$ larger than necessary leads to a large variance of estimation, which is essentially a manifestation of overfitting \cite{chiuso2007role}.

\bibliographystyle{plain}
\bibliography{ref_bib}           



\end{document}